\documentclass[12pt,a4paper]{amsart}
\usepackage[utf8]{inputenc}
\usepackage{amsmath,amssymb,amsthm,xcolor,qtree,bbm}
\usepackage[margin=2.5cm]{geometry}
\usepackage{hyperref}
\title[Banach spaces with Lebesgue property]{Banach spaces with the Lebesgue property of Riemann integrability} 
\author{Harrison Gaebler} 
\author{B\"unyamin Sar\i}
\address{Department of Mathematics, University of North Texas, 1155 Union Circle 311430
Denton, Texas 76203-5017}
\email{harrison.gaebler@unt.edu, bunyamin.sari@unt.edu}

\newcommand{\N}{\mathbb{N}}

\newcommand{\bes}{\begin{equation*}}
\newcommand{\ees}{\end{equation*}}
\newcommand{\ep}{\varepsilon}
\newcommand{\ms}{\min{\text{supp}}}
\newcommand{\Ms}{\max{\text{supp}}}
\newcommand{\supp}{\text{supp}}

\theoremstyle{definition}
\newtheorem{definition}{Definition}[section]
\theoremstyle{plain}
\newtheorem{theorem}[definition]{Theorem}
\newtheorem{lemma}[definition]{Lemma}
\newtheorem{proposition}[definition]{Proposition}

\begin{document}

\maketitle
\begin{abstract}
A Banach space is said to have the \textit{Lebesgue property} if every Riemann-integrable function $f:[0,1]\to X$ is Lebesgue almost everywhere continuous. We give a characterization of the Lebesgue property in terms of a new sequential asymptotic structure that is strictly between the notions of spreading and asymptotic models. We also reproduce an apparently lost theorem of Pelczynski and da Rocha Filho that a subspace $X\subset L_{1}[0,1]$ has the Lebesgue property if every spreading model of $X$ is equivalent to the unit vector basis of $\ell_{1}$.
\end{abstract}
%%%%%%%%%%%%%%%%%%%%%%%%%%%%%%%%%%%%%%%%%%%%%%%%%%%%%%%%%%%%%%%%%%%%%%%%%%%
\section{Introduction}{\label{Sec1}}
Let $X$ be a Banach space. As is the case for scalar-valued functions, a bounded and Lebesgue almost everywhere continuous function $f:[0,1]\to X$ is Riemann-integrable. The converse statement, on the other hand, fails for every classical Banach space except for $X=\ell_{1}$. Indeed, it is easy to see that the function $f:[0,1]\to c_{0}$ (or $\ell_{p}$ for $1<p<\infty$) that maps the rationals to the vectors of the canonical unit vector basis and the irrationals to zero is Riemann-integrable with integral zero. If $X$ does satisfy the condition that every Riemann-integrable function $f:[0,1]\to X$ is Lebesgue almost everywhere continuous, then $X$ is said to have {\em the Lebesgue property.} The problem of determining which classes of Banach spaces have the Lebesgue property is fairly old, dating back at least to the publication in 1927 of \cite{Graves} by Graves that contains what is, as far as we know, the first example that is explicitly noted in print of a Riemann-integrable function $f:[0,1]\to L_{\infty}[0,1]$ that is everywhere discontinuous. However, the literature on this problem is convoluted and somewhat incomplete. For instance, there are several unpublished claims as well as unnoticed solutions to "open" problems. We aim to both clarify and unify in this paper all of the existing literature on the Lebesgue property.

As the above paragraph suggests, Banach spaces with the Lebesgue property have some proximity to $\ell_1$. A result attributed by Haydon in \cite{Haydon} to Pelczynski and da Rocha Filho asserts that {\em every spreading model of $X$ is equivalent to the unit vector basis of $\ell_{1}$  if $X$ has the Lebesgue property}.\footnote{The Lebesgue property is called Darboux property in \cite{Haydon}.} A proof this fact is given by Naralenkov in \cite{Naralenkov}. It is moreover noted in \cite{Haydon} that Pelczynski and da Rocha Filho showed that the converse statement holds if $X\subset L_{1}[0,1]$, however, the original proof of this fact appears to have been lost. Notably, the concept of {\em uniformly separable types} is introduced in \cite{Haydon} and it is proved there that {\em a stable Banach space with uniformly separable types has the Lebesgue property if all spreading models are equivalent to the unit vector basis of $\ell_1$ (equivalently, if $X$ is a Schur space).} While it is not stated in \cite{Haydon}, an unpublished result due to Haydon (see \cite{Chaatit}) is that the space of types on $L_1[0,1]$ is uniformly separable. These two results of Haydon therefore generalize the apparently lost $L_{1}[0,1]$ result of Pelczynski and da Rocha Filho, and it is worth noting also that the separability of the space of types is not an isomorphic invariant (see p. 485 of \cite{Haydon-Maurey})! We give in Section 5 of this paper an alternate and direct proof of the $L_{1}[0,1]$ result of Pelczynski and da Rocha Filho so that it is properly documented in the literature.

In addition, da Rocha Filho shows in \cite{dRF} that both $\ell_{1}$ and the Tsirelson space $T$ have the Lebesgue property. These two results follow from the same argument (see the 1991 survey paper \cite{Gordon} by Gordon as well) and, in particular, this argument is slightly refined in \cite{Naralenkov} to show that a Banach space $X$ that has a basis $(e_{i})$ has the Lebesgue property if it is asymptotic-$\ell_{1}$ with respect to $(e_{i})$. This argument is then further refined by Gaebler in \cite{Gaebler} to show that a Banach space $X$ has the Lebesgue property if it is asymptotic-$\ell_{1}$ in the coordinate-free sense of \cite{MMT}. There is also in \cite{Gaebler} some cursory speculation that the converse statement might hold. However, we show that {\em there is a Banach space $X_{iw}$ with the Lebesgue property which does not contain an asymptotic-$\ell_1$ subspace (Section 4.1).} The space $X_{iw}$ is constructed by Argyros and Motakis, and it has the stronger property that every subspace admits an asymptotic model that is equivalent to the unit vector basis of $c_0$.

It turns out that the Lebesgue property can be characterized in terms of a sequential asymptotic structure. Let $(e_i)$ be a normalized basic sequence in a Banach space $X$. We say that $(e_i)$ is Haar-$\ell^+_1$ sequence if for every Haar  system of partitions $(A^{n}_{j})_{j=0,n\in\N}^{2^{n}-1}$ of $\N$, there exists a constant $C\ge 1$ such that
\bes \frac{1}{C} \leq \lim_{n\to\infty}\sup\left\{\frac{1}{2^{m}}\left\Vert\sum_{j=0}^{2^{m}-1}e_{i_{j}}\right\Vert :\; m\geq n \text{ and } 2^{m}\leq i_{j}\in A^{m}_{j}\right\}. \ees
A system of partitions $(A^{n}_{j})_{j=0,n\in\N}^{2^{n}-1}$ of $\N$ is a Haar system if the ``$n^{\text{th}}$ level" sets  $(A^{n}_{j})_{j=0}^{2^{n}-1}$ form a partition of $\N$ into infinite sets and if the level $n+1$ sets are then obtained by partitioning each of the $n^{\text{th}}$ level sets into two further infinite subsets. Next, we prove (Theorem \ref{MR1}) that {\em a Banach space $X$ has the Lebesgue property if and only if every normalized basic sequence $(e_i)$ in $X$ is a Haar-$\ell^+_1$ sequence.} It follows straightforwardly that every spreading model of $X$ is equivalent to $\ell_{1}$ if all basic sequences are Haar-$\ell^+_1$ and that $X$ has this property if, for every normalized basic sequence $(e_{i})$ in $X$, every asymptotic model of $X$ with respect to an array of normalized block bases of $(e_{i})$ is equivalent to $\ell_{1}$. More interestingly, we show also that the Haar-$\ell^+_1$ condition is in fact strictly between these two spreading and asymptotic model conditions.

After obtaining the above characterization, however, we learned that our result is not new and is essentially the same result that was proved by Pizzotti in her unpublished 1989 PhD thesis (in Portuguese) \cite{Pizzotti}. This seems to have gone unnoticed in the literature, as the later papers such as \cite{Rodriguez} and \cite{KKS} state it as open problem. We thank G. Martinez-Cervantes for sending us the reference \cite{Pizzotti} and we also thank the anonymous referee for pointing out an error in the first draft of this paper.

%%%%%%%%%%%%%%%%%%%%%%%%%%%%%%%%%%%%%%%%%%%%%%%%%%%%%%%%%%%%%%%%%%%%%%%%%%%%%%
\section{Background Material}{\label{Sec2}}
This section contains the background material that is required in order to prove the main results of this paper.
\subsection{Riemann Integration in Banach Spaces}{\label{SubSec2.1}}
The Riemann integral for Banach-valued functions on $[0,1]$ mirrors closely the familiar Riemann integral for scalar-valued functions on $[0,1]$, with the exception that a Banach space $X$ need not have the Lebesgue property. A thorough account of the Riemann integral for Banach-valued functions on $[0,1]$ that includes its comparison to the Darboux integral for Banach-valued functions on $[0,1]$ is found in \cite[Section 2]{Gaebler}. Here, the focus is instead only on the background material that is explicitly required in order to prove the main results of this paper.

Let $f:[0,1]\to X$ be given and recall that a finite and strictly increasing sequence of real numbers $P=(p_{i})_{i=0}^{d}$ is said to be a partition of $[0,1]$ if $p_{0}=0$ and if $p_{d}=1$. The length of the $i^{\text{th}}$ sub-interval of $P$ is $\Delta_{P}(i)=p_{i}-p_{i-1}$ and the number $\pi(P)=\max\{\Delta_{P}(i) \mid 1\leq i\leq d\}$ is then said to be the mesh size of $P$. Next, if $T=(t_{i})_{i=1}^{d}$ is a sequence of real numbers such that $t_{i}\in[p_{i-1},p_{i}]$ (resp. $t_{i}\in(p_{i-1},p_{i})$), then the pair $(P,T)$ is said to be a tagged (resp. interior tagged) partition of $[0,1]$ and the vector
\bes S_{f}(P,T)=\sum_{i=1}^{d}\Delta_{P}(i)f(t_{i}) \ees
is said to be the Riemann sum of $f$ over $(P,T)$. The typical notion of the Riemann integrability of $f$ is now as follows.
\begin{definition}{\label{Rint}}
The function $f:[0,1]\to X$ is said to be Riemann-integrable if there exists a vector $x_{f}\in X$ such that, for every $\ep>0$, there exists a number $\delta=\delta(\ep)>0$ so that $\|x_{f}-S_{f}(P,T)\|<\ep$ for every tagged partition $(P,T)$ of $[0,1]$ with $\pi(P)<\delta$.
\end{definition}
\noindent The set $\mathcal{R}([0,1],X)$ of Riemann-integrable $X$-valued functions on $[0,1]$ is easily seen to be a subspace of the vector space $\mathcal{B}([0,1],X)$ of norm-bounded functions on $[0,1]$ and, moreover, the vector $x_{f}\in X$ is unique for a given $f\in\mathcal{R}([0,1],X)$. In particular, the well-defined linear function $f\mapsto x_{f}$ is said to be the Riemann integral for $X$-valued functions on $[0,1]$. It is worth noting as well that Definition \ref{Rint} is independent of the distinction between tagged and interior tagged partitions. Moreover, as the next proposition shows, one can conveniently take equally-sized dyadic partitions in the definition of Riemann integral.  See \cite{Gaebler} and \cite{Gordon} for other useful equivalent definitions.

\begin{proposition}{\label{prop2.1.2}}
Let $f:[0,1]\to X$ be given. Then, $f\in\mathcal{R}([0,1],X)$ if and only if there exists a vector $x_{f}\in X$ such that, for every $\ep>0$, there is a positive integer $n=n(\ep)\in\N$ so that $\|x_{f}-\frac{1}{2^{m}}\sum_{i=1}^{2^{m}}f(t_{i})\|<\ep$ for every $m\geq n$ and for every choice $t_{i}\in\left(\frac{i-1}{2^{m}},\frac{i}{2^{m}}\right)$ of interor tags.
\end{proposition}
\begin{proof}
The ``if" direction is trivial so it remains only to prove the ``only if" direction. Let $\ep>0$ be given and choose $n=n(\ep)\in\N$ so that $\|\frac{1}{2^{m}}\sum_{i=1}^{2^{m}}f(t_{i})\|<\ep$ for every $m\geq n$ and for every choice $t_{i}\in\left(\frac{i-1}{2^{m}},\frac{i}{2^{m}}\right)$ of interior tags, where we have assumed without loss of generality (by translating) that $x_{f}=0$. 

Note that $f$ is bounded and put $M=\sup_{t\in[0,1]}\|f(t)\|$. Then, let $P=(p_{j})_{j=0}^{d}$ be a partition of $[0,1]$ such that $\pi(P)<\frac{\ep}{M2^{n+1}}$, define the set $B=\left\{j \mid \frac{i}{2^{n}}\in [p_{j-1},p_{j}] \text{ for some } i \right\}$ so that $|B| \leq 2^{n+1}-2$, and let $(t_{i})_{i=1}^{2^{n}}$ and $(s_{j})_{j=1}^{d}$ be interior tags. It now follows that
\begin{align*} &\left\Vert \frac{1}{2^{n}}\sum_{i=1}^{2^{n}}f(t_{i}) - \sum_{j=1}^{d}\Delta_{P}(j)f(s_{j})\right\Vert \leq \left\Vert\sum_{i=1}^{2^{n}}\left(\frac{f(t_{i})}{2^{n}}-\sum_{j\in A_{i}}\Delta_{P}(j)f(s_{j})\right)\right\Vert + \left\Vert\sum_{j\in B}\Delta_{P}(j)f(s_{j})\right\Vert \\ &\qquad \leq \left\Vert \sum_{i=1}^{2^{n}}\left(\sum_{j\in A_{i}}\Delta_{P}(j)+\rho(i)\right)f(t_{i}) -\sum_{j\in A_{i}}\Delta_{P}(j)f(s_{j})\right\Vert + \ep  \\ &\qquad \leq \left\Vert\sum_{i=1}^{2^{n}}\sum_{j\in A_{i}}\Delta_{P}(j)(f(t_{i})-f(s_{j})) \right\Vert + \left\Vert\sum_{i=1}^{2^{n}}\rho(i)f(t_{i})\right\Vert + \ep \\ &\qquad \leq \left\Vert\sum_{i=1}^{2^{n}}\sum_{j\in A_{i}} 2^{n}\Delta_{P}(j)\frac{1}{2^{n}}(f(t_{i})-f(s_{j}))\right\Vert + 2\ep \qquad (\#) \end{align*}
where $A_{i}=\left\{j \mid [p_{j-1},p_{j}]\subset\left(\frac{i-1}{2^{n}},\frac{i}{2^{n}}\right)\right\}$ and $\frac{1}{2^{n}}=\sum_{j\in A_{i}}\Delta_{P}(j)+\rho(i)$ with $0\leq\rho(i)\leq 2\pi(P)$. We may assume by a small perturbation of the coefficients $2^{n}\Delta_{P}(j)$ that
\bes \sum_{i=1}^{2^{n}}\sum_{j\in A_{i}}2^{n}\Delta_{P}(j)\frac{1}{2^{n}}(f(t_{i})-f(s_{j})) \in\sum_{i=1}^{2^{n}}\text{co}(Y_{i})=\text{co}\left(\sum_{i=1}^{2^{n}}Y_{i}\right)\ees
where $Y_{i}=\left\{\frac{1}{2^{n}}(f(t)-f(t')) \;\middle\vert\; t,t'\in\left(\frac{i-1}{2^{n}},\frac{i}{2^{n}}\right)\right\}$ and where $\text{co}(\cdot)$ is the convex hull operator. Then, $\|x\|<2\ep$ by hypothesis if $x\in\text{co}\left(\sum_{i=1}^{2^{n}}Y_{i}\right)$ so that $(\#)\leq 4\ep$ and the proof is complete.
\end{proof}

Recall that a sequence $(e_{i})$ in $X$ is said to be basic if, for every $x\in\overline{\text{span}\{e_{i}\}}:=[e_{i}]$, there exists a unique sequence $(a_{i})$ of scalars such that $x=\sum_{i}a_{i}e_{i}$ or, equivalently, if there exists a constant $K_{b}\geq 1$ such that $\left\Vert\sum_{i=1}^{n}a_{i}e_{i}\right\Vert\leq K_{b}\left\Vert\sum_{i=1}^{m}a_{i}e_{i}\right\Vert$ for every $m\geq n$ and for every sequence $(a_{i})$ of scalars. In particular, a basic sequence $(e_{i})$ is said to be a basis for $X$ if $X=[e_{i}]$. The set $\text{supp}(x)=\{i\in\N \mid a_{i}\neq 0\}$ is then said to be the support of a vector $x\in[e_{i}]$ and a sequence $(x_{j})$ of non-zero vectors in $[e_{i}]$ is said to be a block basis of $(e_{i})$ if, for every $j\in\N$, $x_{j}$ has finite support and $\max\text{supp}(x_{j})<\min\text{supp}(x_{j+1})$.

The second result we include in this subsection is the following proposition due to Haydon and Odell which says that the failure of the Lebesgue property in a space $X$ is witnessed by a `Dirichlet's function' which vanishes outside the set of dyadic rationals. As observed by Pizzotti \cite[Thm. 1.14]{Pizzotti} the range of the function can also be chosen to be a normalized basic sequence in the space. 

\begin{proposition}{\label{Haydon-Odell}}
Assume that $X$ does not have the Lebesgue property. Then, there exists a normalized basic sequence $(x_{j})$ in $X$ so that the function $\mathcal{F}:[0,1]\to X$ that is defined by
\bes \mathcal{F}(t)=\begin{cases} x_{2^{n}+i} &\mbox{if } n\in\N\cup\{0\},\text{ }0\leq i\leq 2^{n}-1,\text{ and } t=\frac{2i+1}{2^{n+1}} \\ 0 &\text{otherwise} \end{cases} \ees
is Riemann-integrable. Moreover, $(x_{j})$ can be taken to be a normalized block basis of $(e_{i})$ if $(e_{i})$ is a basis for $X$.
\end{proposition}

This result 
was unpublished, but mentioned by Haydon in \cite{Haydon} and a proof was given by Pizzotti, \cite[Thm. 1.14]{Pizzotti}. We will reproduce a proof below for completeness.

First, we have the following auxiliary proposition.

\begin{proposition}\label{dirichlet}
Assume that $X$ does not have the Lebesgue property. Then, there exists a countable and dense subset $\mathcal{A}\subset[0,1]$ and a Riemann-integrable function $f:[0,1]\to X$ such that $\|f(t)\|=1$ if $t\in \mathcal{A}$ and $f(t)=0$ if $t\notin \mathcal{A}$.
\end{proposition}

\begin{proof} The desired function will be defined in few steps. 
Since $X$ does not have the Lebesgue property, there exists a Riemann-integrable function $f_1:[0,1]\to X$ that is discontinuous on a set of $F$ of positive Lebesgue measure. Then, define for each $n\in\N$ the Lebesgue measurable set
\bes F_{n}=\left\{ t\in[0,1] \;\middle\vert\; \inf_{\delta>0}\sup\left\{\|f_1(s)-f_1(s')\| \mid s,s'\in[t-\delta,t+\delta]\cap[0,1]\right\}\geq\frac{1}{n}\right\} \ees
and note that $F=\bigcup_{n=1}^{\infty}F_{n}$ so that $\mu(F_{n_{0}})>0$ for some $n_{0}\in\N$ where $\mu$ is the Lebesgue measure. Now, write $I^{j}_{k}=\left(\frac{j}{2^{k}},\frac{j+1}{2^{k}}\right)$ for each $k\in\N$ and for each $0\leq j\leq 2^{k}-1$ and note that, if $I^{j}_{k}\cap F_{n_{0}}\neq\emptyset$, then there exist points $x_{j,k},y_{j,k}\in I^{j}_{k}$ such that $\|f_1(x_{j,k})-f_1(y_{j,k})\|\geq\frac{1}{2n_{0}}$. For each such $I^{j}_{k}$ fix $x_{j,k},y_{j,k}\in I^{j}_{k}$ as above. Let $m_{j,k}$ be the midpoint of $I^{j}_{k}$ and define $f_2:[0,1]\to X $ by
\bes f_2(t)=\begin{cases} f_1(x_{j,k})-f_1(y_{j,k}) &\mbox{if } t=m_{j,k} \text{ where } I^{j}_{k}\cap F_{n_{0}}\neq\emptyset \\ 0 &\text{otherwise} \end{cases} \ees
It follows from Proposition \ref{prop2.1.2} that $f_2$ is Riemann-integrable because $f_1$ is Riemann-integrable. Note that $f_2$ is discontinuous on $F_{n_{0}}$. Let $A=\{t\in[0,1]: f_2(t)\neq 0\}$. Note $A\subseteq D$, where $D$ is the set of dyadic rationals in $(0,1)$, and in particular, $f_2$ vanishes on a dense subset of $[0,1]$

Next, we replace the function $f_2$ with one whose support is countable and {\em dense} in $[0,1]$. Let $U$ be the union the dyadic intervals in $[0,1]$ where $f_2$ vanishes and note that $\mu([0,1]\setminus U)>0$ because $f_2$ is discontinuous on a set of positive Lebesgue measure. Now, define the function $\alpha:[0,1]\to [0,1]$ by
\bes \alpha(t)=\frac{\mu\left([0,t]\cap\left([0,1]\setminus U\right)\right)}{\mu([0,1]\setminus U)} \ees
and note that $\alpha$ is continuous on $[0,1]$, and $\alpha(0)=0$ and $\alpha(1)=1$. It therefore follows from the intermediate value theorem that, for each $x\in[0,1]$, there exists some $t_{x}\in[0,1]$ so that $x=\alpha(t_{x})$. Then, define the function $f_3:[0,1]\to X$ by $f_3(x)=f_2(t_{x})$ and note that $f_3(x)=0$ if $x\notin\alpha(A)$. Note also that $\mathcal{A}:=\alpha(A)$ is countable because $A$ is countable. We will prove that $f_3$ is Riemann-integrable and that $\mathcal{A}$ is dense in $[0,1]$.

Recall that $f_2$ is Riemann-integrable. Let $\ep>0$ be given and choose $n\in\N$ so that there is the `sign-unconditional' estimate
\bes  \left\Vert \frac{1}{2^{m}}\sum_{i=1}^{2^{m}}\eta_{i}(f_2(u_{i})-f_2(v_{i}))\right\Vert \leq \frac{\mu([0,1]\setminus U)}{2}\ep.
\ees
for every $m\geq n$, for every choice of signs $\eta_{i}\in\{-1,1\}$, and for every choice of tags $u_{i},v_{i}\in\left[\frac{i-1}{2^{m}},\frac{i}{2^{m}}\right]$. Note that 
\bes \left[\alpha\left(\frac{i}{2^{m}}\right)-\alpha\left(\frac{i-1}{2^{m}}\right)\right]\mu([0,1]\setminus U)=\mu\left(\left(\frac{i-1}{2^{m}},\frac{i}{2^{m}}\right)\cap([0,1]\setminus U)\right)\leq\frac{1}{2^{m}}. \ees
Thus from the sign-unconditionality of the above expression it follows that
\bes \left\Vert \sum_{i=1}^{2^{m}}(f_2(u_{i})-f_2(v_{i}))\left[\alpha\left(\frac{i}{2^{m}}\right)-\alpha\left(\frac{i-1}{2^{m}}\right)\right]\mu([0,1]\setminus U)\right\Vert \leq \mu([0,1]\setminus U)\ep \ees
and thus,
\bes \left\Vert \sum_{i=1}^{2^{m}}(f_2(u_{i})-f_2(v_{i}))\left[\alpha\left(\frac{i}{2^{m}}\right)-\alpha\left(\frac{i-1}{2^{m}}\right)\right]\right\Vert \leq \ep \ees
for every choice of tags $u_{i},v_{i}\in\left[\frac{i-1}{2^{m}},\frac{i}{2^{m}}\right]$ after dividing by $\mu([0,1]\setminus U)$. Let $\tau_{i}=\frac{i}{2^{m}}$ for each $0\leq i\leq 2^{m}$ and let $0=\tau_{i_{1}}<\ldots<\tau_{i_{r}}=1$ be the at most $2^{m}+1$ values such that, for each $1\leq k\leq r$, we have that $\alpha(\tau_{i_{k-1}})<\alpha(\tau_{i_{k}})$ and $\alpha(\tau_{i})=\alpha(\tau_{i_{k-1}})$ for each $\tau_{i}\in[\tau_{i_{k-1}},\tau_{i_{k}})$. Then, $(\alpha(\tau_{i_{k}}))_{k=1}^{r}$ is a partition of $[0,1]$ and, for every choice of tags $x_{k},y_{k}\in[\alpha(\tau_{i_{k-1}}),\alpha(\tau_{i_{k}})]=[\alpha(\tau_{i_{k}-1}),\alpha(\tau_{i_{k}})]$, there exist points $u_{k},v_{k}\in[\tau_{i_{k}-1},\tau_{i_{k}}]$ so that $x_{k}=\alpha(u_{k})$ and $y_{k}=\alpha(v_{k})$. It therefore follows that
\bes \left\Vert \sum_{k=1}^{r}(f_3(x_{k})-f_3(y_{k}))(\alpha(\tau_{i_{k}})-\alpha(\tau_{i_{k-1}}))\right\Vert = \left\Vert \sum_{k=1}^{r}(f_2(u_{k})-f_2(v_{k}))(\alpha(\tau_{i_{k}})-\alpha(\tau_{i_{k}-1}))\right\Vert\leq \ep \ees
for every choice of tags $x_{k},y_{k}\in[\alpha(\tau_{i_{k-1}}),\alpha(\tau_{i_{k}})]$, and this proves that $f_3$ is Riemann-integrable (e.g. by \cite[Thm. 2.2.3.]{Gaebler}. It remains to show that $\mathcal{A}$ is dense in $[0,1]$.

Let $t\in([0,1]\setminus U)$ and let $V$ be an arbitrary open set that contains $t$. Then, there is some dyadic interval $I$ so that $t\in I\subset V$ and, by construction, we have that $A\cap I\neq\emptyset$. It follows that $A\cap V\neq\emptyset$ and therefore $t\in \overline{A}$ because $V$ is arbitrary.

Next, let $J=\{x\in[0,1] \mid \alpha^{-1}(\{x\}) \text{ is not a singleton}\}$ and note that $\alpha(U)\subset J$. Note also that $\alpha^{-1}(\{x\})$ is a closed interval with positive Lebesgue measure for each $x\in J$ because $\alpha$ is continuous and is (not necessarily strictly) increasing. It follows that $J$ is countable because $[0,1]$ is the union of at most countably many closed intervals with positive Lebesgue measure.

Now, if $x\notin J$, then there exists a point $t_{x}\in[0,1]\setminus U$ so that $x=\alpha(t_{x})$. Thus, $t_{x}\in\overline{A}$ so that $x\in\alpha(\overline{A})\subset\overline{\alpha(A)}$ where the inclusion follows from the continuity of $\alpha$. We have therefore shown that $[0,1]\setminus J\subset\overline{\mathcal{A}}$ so that $\mu(\overline{\mathcal{A}})=1$ and it follows that $\mathcal{A}$ is dense in $[0,1]$.

Finally, define 
\bes f(t)=\begin{cases} \frac{f_3(t)}{\|f_3(t)\|} &\mbox{if } t\in\mathcal{A}  \\ 0 &\text{otherwise} \end{cases} \ees

It follows from the fact that $f_3$ is Riemann integrable and it vanishes on a dense set that $f$ is Riemann integrable. Indeed, let $s=1+\sup_{t\in[0,1]}\|f_3(t)\|$. Then, let $\ep>0$ be given and and choose $n\in\N$ so that $\left\Vert\frac{1}{2^{m}}\sum_{i=1}^{2^{m}}(f_3(t_{i})-f_3(t'_{i}))\right\Vert<\frac{\ep}{2s}$ for every $m\geq n$ and for every choice of interior tags $t_{i},t'_{i}\in\left(\frac{i-1}{2^{m}},\frac{i}{2^{m}}\right)$. It therefore follows that $\left\Vert\frac{1}{2^{m}}\sum_{i=1}^{2^{m}}\ep_{i}f_3(t_{i})\right\Vert<\frac{\ep}{2s}$ for every $m\geq n$, for every choice of interior tag $t_{i}\in\left(\frac{i-1}{2^{m}},\frac{i}{2^{m}}\right)$, and for every choice of sign $\ep_{i}\in\{-1,1\}$ because $f_3$ vanishes on a dense subset of $[0,1]$. Now, again by the sign-unconditionality of the above expression, we have
\bes \left\Vert \frac{1}{2^{m}}\sum_{i=1}^{2^{m}}\frac{f_3(t_{i})}{\|f_3(t_{i})\|}\right\Vert  \leq s\cdot 2\frac{\ep}{2s}=\ep  \ees
 so that the proof is complete.
\end{proof}

Note that the Lebesgue property is {\em separably-determined} by the above proposition.

Now, we are ready to prove the Proposition \ref{Haydon-Odell} using the above result and a standard gliding hump argument.

\begin{proof}[Proof of Proposition \ref{Haydon-Odell}]
Assume without loss of generality that $X$ is separable and regard $X$ as a subspace of $C[0,1]$. Let $(e_{i})$ be the monotone Schauder basis of $C[0,1]$. By Proposition \ref{dirichlet} there exists a countable and dense subset $\mathcal{A}\subset[0,1]$ and a Riemann-integrable function $f:[0,1]\to X$ such that $\|f(t)\|=1$ if $t\in\mathcal{A}$ and $f(t)=0$ if $t\notin\mathcal{A}$.

Let $P_{n}$ be the basis projection onto the span of the first $n$ basis vectors. Then, since $P_{n}$ is continuous, the composition $P_{n}\circ f$ is Riemann integrable, and thus Lebesgue almost-everywhere continuous since its range is finite-dimensional. Observe that for all $\ep>0$ and for every $n\in\N$ and for every interval $(a,b)\subset[0,1]$, there exists $t\in(a,b)\setminus\mathcal{A}$ and $s\in(a,b)\cap\mathcal{A}$ such that $\|(P_{n}\circ f)(t)-(P_{n}\circ f)(s)\|=\|(P_{n}\circ f)(s)\|<\ep$ from the Lebesgue almost-everywhere continuity of $P_{n}\circ f$ and the density in $[0,1]$ of $\mathcal{A}$, whereas $\|f(s)\|=1$.

Let  $I_{1}=(0,1)$, $I_{2}=(0,\frac{1}{2})$, $I_{3}=(\frac{1}{2},1)$, and so on, and let $d_{1}=\frac{1}{2}$, $d_{2}=\frac{1}{4}$, $d_{3}=\frac{3}{4}$, and so on be the dyadic rationals so that $I_{2^{n}+i}=\left(\frac{i}{2^{n}},\frac{i+1}{2^{n}}\right)$ and $d_{2^{n}+i}=\frac{2i+1}{2^{n+1}}$ for each $n\in\N\cup\{0\}$ and for each $0\leq i\leq 2^{n}-1$. We will define a new function $g:[0,1]\to C[0,1]$ so that $g$ vanishes outside of $D$. Let $s_{1}\in I_{1}$ so that $\|f(s_{1})\|=1$  by the density of $\mathcal{A}$ and choose $n_{1}\in\N$ so that $\|f(s_{1})-P_{n_{1}}(f(s_{1}))\|<\frac{1}{4}$. Then, define $g(d_{1})=P_{n_{1}}(f(s_{1}))$ and note that 
\bes \|g(d_{1})\|\geq \|f(s_{1})\|-\frac{1}{4}>\frac{1}{2} \text{ while } \|g(d_{1})-f(s_{1})\|<\frac{1}{4}.\ees

Next, choose $s_{2}\in I_{2}$ so that $\|P_{n_{1}}(f(s_{2}))\|<\frac{1}{16}$ and $\|f(s_{2})\|=1$ by the density of $\mathcal{A}$, and choose $n_{2}\in\N$ with $n_{2}>n_{1}$ so that $\|f(s_{2})-P_{n_{2}}(f(s_{2}))\|<\frac{1}{16}$. Define $g(d_{2})=(P_{n_{2}}-P_{n_{1}})(f(s_{2}))$ and note that $\|g(d_{2})-f(s_{2})\|\leq \frac{1}{8}$ so that
\bes \|g(d_{2})\| \geq \|f(s_{2})\|-\|f(s_{2})-g(d_{2})\| \geq 1-\frac{1}{8}>\frac{1}{2} \ees
and continue this gliding hump process to extract an increasing sequence $(n_{j})$ of positive integers and a sequence $(s_{j})$ of points with $s_{j}\in I_{j}$ so that $\|g(d_{j})\|>\frac{1}{2}$ and $\|f(s_{j})-g(d_{j})\|<\frac{1}{2^{j+1}}$ for each $j\in\N$. Let $g(t)=0$ if $t\neq d_{j}$ for any $j\in\N$.

Now, define $\mathcal{F}:[0,1]\to X$ by $\mathcal{F}(t)=f(s_{j})$ if $t=d_{j}$ and $\mathcal{F}(t)=0$ otherwise. It then follows from Proposition \ref{prop2.1.2} that $\mathcal{F}$ is Riemann-integrable because $d_{j}$ is the midpoint of the dyadic interval $I_{j}$ from which $s_{j}$ is chosen and because $f$ is Riemann-integrable. Moreover, $(f(s_j))_j$ is a normalized basic sequence in $X$ since, by construction and the principle of small perturbations, it is equivalent to the block basis $(g(d_j))_j$ in $C[0,1]$.

To see the moreover statement, suppose $(e_i)$ is a basis for $X$. In the above argument, take $\mathcal{F}:[0,1]\to X$ by $\mathcal{F}(t)=\frac{g(d_{j})}{\|g(d_{j}\|}$ if $t=d_{j}$ and $\mathcal{F}(t)=0$ otherwise. Again $\mathcal{F}$ is Riemann-integrable since $(g(d_j))_j$ is equivalent to $(f(s_j))_j$.
\end{proof}

\subsection{Haar systems of partitions of $\N$}{\label{SubSec2.2}}
Our characterization of the Lebesgue property is based on partitions of $\N$ that mimic the structure of the dyadic rationals in $[0,1)$. Note first that the set $D'=D\cup\{0\}$ of dyadic rational numbers in $[0,1)$ is equal to $\bigcup_{n=1}^{\infty}\left\{\frac{j}{2^{n}} \mid 0\leq j \leq 2^{n}-1\right\}$ and consider the dyadic tree $\mathcal{D}$ that is given by 
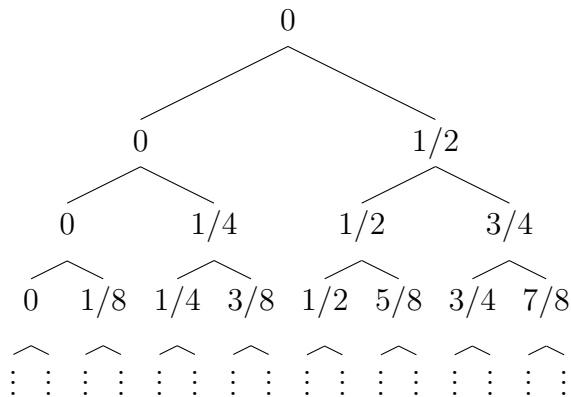
\bes \Tree [.0 [ [ [{\vdots} {\vdots} ].0 [{\vdots} {\vdots} ].1/8 ].0 [[{\vdots} {\vdots} ].1/4 [{\vdots} {\vdots} ].3/8 ].1/4 ].0  [[[{\vdots} {\vdots} ].1/2 [{\vdots} {\vdots} ].5/8 ].1/2 [[{\vdots} {\vdots} ].3/4 [{\vdots} {\vdots} ].7/8 ].3/4 ].1/2 ] \ees
and that collects the members of $D'$ according to their height $n\in\N$. In particular, for every $d^{n}_{j}=\frac{j}{2^{n}}\in D'$ (we assume that $j$ is odd so that $\frac{j}{2^{n}}$ is in lowest terms), there exists a unique branch in $\mathcal{D}$ such that the $n^{th}$ node of this branch is $d^{n}_{j}$. We now re-create this tree structure for certain infinite subsets of $\N$. 
\begin{definition}{\label{dyadic_system}}
A collection $(A^{n}_{j})_{j=0,n\in\N}^{2^{n}-1}$ of infinite subsets of $\N$ is said to be a {\em Haar system of partitions} of $\N$ if:
\begin{enumerate}
\item for every $n\in\N$,  $\bigcup_{j=0}^{2^{n}-1}A^{n}_{j}=\N$ and $A^{n}_{j}\cap A^{n}_{j'}=\emptyset$ if $j\neq j'$
\item for every $n\in\N$ and $0\le j<2^n-1$, $A^{n}_{j}=A^{n+1}_{2j}\cup A^{n+1}_{2j+1}$
\end{enumerate}
 where the second requirement is simply the assertion that the level $n+1$ subsets are given by partitioning every level $n$ subset into two further infinite subsets.
\end{definition}

\noindent A Haar system of partitions $(A^{n}_{j})_{j=0,n\in\N}^{2^{n}-1}$ then can be arranged in a tree that is similar to $\mathcal{D}$ as follows,
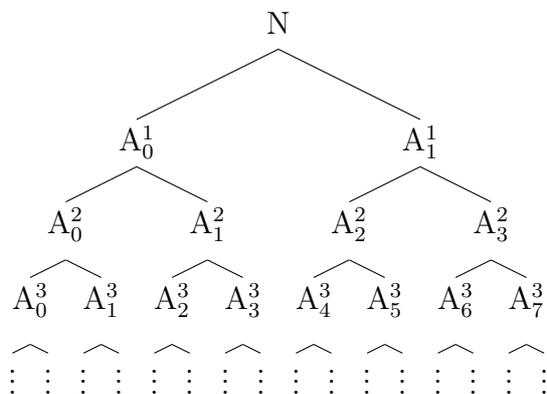
\bes \Tree [.N [ [ [{\vdots} {\vdots} ].A^{3}_{0} [{\vdots} {\vdots} ].A^{3}_{1} ].A^{2}_{0} [[{\vdots} {\vdots} ].A^{3}_{2} [{\vdots} {\vdots} ].A^{3}_{3} ].A^{2}_{1} ].A^{1}_{0}  [[[{\vdots} {\vdots} ].A^{3}_{4} [{\vdots} {\vdots} ].A^{3}_{5} ].A^{2}_{2} [[{\vdots} {\vdots} ].A^{3}_{6} [{\vdots} {\vdots} ].A^{3}_{7} ].A^{2}_{3} ].A_{1}^{1} ] \ees
and we define the map $\sigma$ to be the function that takes $d^{n}_{j}$ to the $(2^{n})^{\text{th}}$ member of $A^{n}_{j}$ (note that $\sigma$ is well-defined because $d^{n}_{j}$ is in lowest terms). Thus $\sigma$ embeds the dyadic rationals in $[0,1)$ into the Haar system of partitions of $\N$ so that, for every $d^{n}_{j}=\frac{j}{2^{n}}$, the unique branch in $\mathcal{D}$ whose $n^{\text{th}}$ node is $d^{n}_{j}$ corresponds to a unique branch in the above tree.  Most importantly, it follows that the dyadic rationals in a given interval of the form $\left[\frac{j-1}{2^{n}},\frac{j}{2^{n}}\right)$ are mapped by $\sigma$ to ``deeper level" subsets under $A^{n}_{j-1}$ on the same bush. This structural observation is key to our characterization of the Lebesgue property.

\begin{definition}{\label{Haar-1}}
We say that a normalized basic sequence $(e_i)$ in a Banach space $X$ is Haar-$\ell^+_1$ sequence if for every Haar  system of partitions $(A^{n}_{j})_{j=0,n\in\N}^{2^{n}-1}$ of $\N$, there exists a constant $C\ge 1$ such that
\bes \frac{1}{C} \leq \lim_{n\to\infty}\sup\left\{\frac{1}{2^{m}}\left\Vert\sum_{j=0}^{2^{m}-1}e_{i_{j}}\right\Vert :\; m\geq n \text{ and } 2^{m}\leq i_{j}\in A^{m}_{j}\right\}. \ees
\end{definition}

\section{A Characterization of Banach spaces with the Lebesgue property}{\label{SubSec3.1}}
We begin with the main theorem of this subsection.
\begin{theorem}{\label{MR1}}
A Banach space $X$ has the Lebesgue property if and only if every normalized basic sequence $(e_i)$ in $X$ is a Haar-$\ell^+_1$ sequence.
\end{theorem}
\begin{proof}
Assume $X$ does not have the Lebesgue property. Then, by Proposition \ref{Haydon-Odell}, there is a normalized basic sequence $(x_{j})$ so that the function $f:[0,1]\to X$ that is defined by
\bes f(t)=\begin{cases} x_{j} &\mbox{if } t=d_{j}\in D \\ 0 &\mbox{if } t\notin D \end{cases} \ees
is Riemann-integrable, where $D$ is as usual the set of dyadic rational numbers in $(0,1)$. Note also that the Riemann integral itself of the Riemann-integrable function $f$ is in this situation equal to zero. Let $(A^{n}_{k})_{k=0, n\in\N}^{2^{n}-1}$ be the Haar system of partitions of $\N$ so that, for every $n\in\N$,  $A^{n}_{0}=\left\{j \;\middle\vert\; d_{j}\in\left(0,\frac{1}{2^{n}}\right)\right\}$, $A^{n}_{1}=\left\{j \;\middle\vert\; d_{j}\in\left[\frac{1}{2^{n}},\frac{2}{2^{n}}\right)\right\},\ldots,$ $A^{n}_{2^{n}-1}=\left\{j \;\middle\vert\; d_{j}\in\left[\frac{2^{n}-1}{2^{n}},1\right)\right\}$. Let $C\ge 1$ be arbitrary and pick $\ep\in(0,\frac{1}{2C})$. Then, by integrability, there exists a positive integer $n=n(\ep)\in\N$ such that $\frac{1}{2^{m}}\left\Vert\sum_{i=1}^{2^{m}}f(t_{i})\right\Vert<\ep$ for every $m\geq n$ and for {\em every choice} $t_{i}\in\left[\frac{i-1}{2^{m}},\frac{i}{2^{m}}\right)$ of tags. It therefore follows that
\bes \frac{1}{2^{m}}\left\Vert\sum_{k=0}^{2^{m}-1}x_{j_{k}}\right\Vert=\frac{1}{2^{m}}\left\Vert\sum_{k=0}^{2^{m}-1}f(d_{j_{k}})\right\Vert<\ep<\frac1C \ees
for $m\geq n$ and {\em for all} $j_{k}\in A^{m}_{k}$, $0\le k\le 2^m-1$. In other words, $(x_j)$ is not a Haar-$\ell^+_1$ sequence.

Conversely, assume that $X$ has the Lebesgue property and let $(e_{i})$ be a normalized basic sequence in $X$. Fix also a Haar system of partitions $(A^{n}_{j})_{j=0,n\in\N}^{2^{n}-1}$ of $\N$. Let $(d_{k})$ be the natural ordering of the dyadic rationals $D$ in $(0,1)$ (that is, 1/2, 1/4, 3/4, 1/8, and so on) and define a function $f:[0,1]\to X$ by 
\bes f(t)=\begin{cases} e_{\sigma(d_{k})} &\mbox{if } t=d_{k}\in D \\ 0 &\mbox{if } t\notin D \end{cases} \ees
where $\sigma$ is the map that is defined in Subsection \ref{SubSec2.2}. It follows that $f$ is discontinuous everywhere and is therefore not Riemann-integrable. This implies, by Proposition \ref{prop2.1.2}, that there exists a constant $c_{f}>0$ such that, for every $n'\in\N$, there is $n\ge n'$ such that  $c_{f}\leq\frac{1}{2^{n}}\left\Vert\sum_{j=1}^{2^{n}}f(t_{j})\right\Vert$ for some interior tags $t_{j}\in\left(\frac{j-1}{2^{n}},\frac{j}{2^{n}}\right)$. For such $n\in\N$ consider the ``level $n$" sets $A^{n}_{0},\ldots, A^{n}_{2^{n}-1}$ of the Haar system of partitions of $\N$, and let $t_{j}\in\left(\frac{j-1}{2^{n}},\frac{j}{2^{n}}\right)$ be interior tags so that $c_{f}\leq\frac{1}{2^{n}}\left\Vert\sum_{j=1}^{2^{n}}f(t_{j})\right\Vert$. Define the set $B=\{j \mid t_{j}=d_{k_{j}}\in D \text{ for some }k_{j}\}$, that is $B=\{j \mid f(t_j)\neq 0\}$. Note that for every $j\in B$, $\sigma(d_{k_{j}})=\sigma\left(\frac{l}{2^{m}}\right)\in A^{m}_{l}\subset A^{n}_{j}$ for some $m>n$ and for some $1\leq l\leq 2^{m}-1$. In particular, $\sigma(d_{k_{j}})$ is by the definition of $\sigma$ at least the $(2^{n})^{\text{th}}$ member of $A^{n}_{j}$. For $j\notin B$ choose any dyadic rational $d_{j}\in\left(\frac{j-1}{2^{n}},\frac{j}{2^{n}}\right)$ such that $\max_{j\in B}\sigma(d_{k_{j}})<\sigma(d_{j})\in A^{n}_{j}$. Then \begin{align*} \frac{1}{2^{n}}\left\Vert\sum_{j\in B}e_{\sigma(d_{k_{j}})}+\sum_{j\notin B}e_{\sigma(d_{j})} \right\Vert &\geq\frac{1}{K_{b}2^{n}}\left\Vert\sum_{j\in B}e_{\sigma(d_{k_{j}})}\right\Vert \\ &= \frac{1}{K_{b}2^{n}}\left\Vert\sum_{j=1}^{2^{n}}f(t_{j})\right\Vert\geq \frac{c_{f}}{K_{b}}\end{align*}
where $K_{b}\geq 1$ is the basis constant for $(e_{i})$ and this implies that $(e_i)$ is a Haar-$\ell^+_1$ sequence with constant $\frac{K_b}{c_f}$, so that the proof is complete.
\end{proof}

\section{Asymptotic structures and Haar-$\ell^+
_1$ sequences}
The Haar-$\ell^+_1$ condition is a sequential asymptotic property. In this section, we explore the connections between other sequential asymptotic structures. We will show that the condition that every basic sequence in $X$ is Haar-$\ell^+_1$ is strictly between the property that every spreading model of $X$ is equivalent to $\ell_{1}$ and the property that every asymptotic model of $X$ with respect to normalized block bases of a given basic sequence is equivalent to $\ell_{1}$. Recall that two basic sequences $(e_{i})$ and $(e'_{i})$ are said to be equivalent if there exists a constant $C>0$ such that, for every $n\in\N$, there is
\bes \frac{1}{C}\left\Vert\sum_{i=1}^{n}a_{i}e'_{i}\right\Vert\leq\left\Vert\sum_{i=1}^{n}a_{i}e_{i}\right\Vert\leq C\left\Vert\sum_{i=1}^{n}a_{i}e'_{i}\right\Vert \ees
for every sequence $(a_{i})$ of scalars. The definition of a spreading model is now given below and it is due to Brunel and Sucheston in \cite{BS}.
\begin{definition}{\label{sm}}
Let $(e_{i})$ be a normalized basic sequence in $X$ and let $(v_{i})$ be a normalized basis for a Banach space $(V,\|\cdot\|_{V})$. If there exist positive real numbers $\delta_{n}\searrow 0$ such that
\bes \left\vert \left\Vert\sum_{j=1}^{n}a_{j}e_{i_{j}}\right\Vert-\left\Vert\sum_{i=1}^{n}a_{i}v_{i}\right\Vert_{V}\right\vert<\delta_{n} \ees
for every choice $n\leq i_{1}<\ldots <i_{n}$ of positive integers and for every sequence $(a_{i})_{i=1}^{n}$ of scalars with $|a_{i}|\leq 1$, then $(v_{i})$ is said to be {\em a spreading model of} $X$ generated by $(e_{i})$.
\end{definition}
\noindent A normalized basic sequence $(e_{i})$ is said to be good if it generates a spreading model $(v_{i})$ of $X$ and, as Brunel and Sucheston show in \cite{BS}, it is a straightforward consequence of Ramsey's theorem that every normalized basic sequence has a subsequence that is good. In addition, a spreading model $(v_{i})$ of $X$ is necessarily $1$-spreading, that is, every subsequence of $(v_{i})$ is $1$-equivalent to $(v_{i})$. The notion of a spreading model is generalized by Halbeisen and Odell \cite{HO} in terms of arrays of basic sequences. An array $(e_{i}^{j})$ of vectors in $X$ is said to be $C$-basic if, for every $j\in\N$, $(e^{j}_{i})$ is a normalized $C$-basic sequence and if, for every choice $n\leq i_{1}< i_{2}<\ldots <i_n$ of positive integers, the diagonal sequence $(e^{j}_{i_{j}})_{j=1}^n$ is $C$-basic. 
\begin{definition}{\label{am}}
Let $(e_{i}^{j})$ be a $C$-basic array in $X$ and let $(v_{i})$ be a normalized basis for a Banach space $(V,\|\cdot\|_{V})$. If there exist positive real numbers $\delta_{n}\searrow 0$ such that
\bes \left\vert \left\Vert \sum_{j=1}^{n}a_{j}e^{j}_{i_{j}}\right\Vert - \left\Vert\sum_{i=1}^{n}a_{i}v_{i}\right\Vert_{V}\right\vert < \delta_{n} \ees
for every choice $n\leq i_{1}<\ldots<i_{n}$ of positive integers and for every sequence $(a_{i})_{i=1}^{n}$ of scalars with $|a_{i}|\leq 1$, then $(v_{i})$ is said to be {\em an asymptotic model of} $X$ generated by $(e_{i}^{j})$.
\end{definition}
\noindent A $C$-basic array $(e_{i}^{j})$ is said to be good if it generates an asymptotic model $(v_{i})$ of $X$ and, as Halbeisen and Odell show in \cite{HO}, it is again a straightforward consequence of Ramsey's theorem that every $C$-basic array has a $C$-basic subarray that is good. A subarray is obtained by passing to subsequences in each row. Unlike a spreading model, an asymptotic model $(v_{i})$ need not be $C$-spreading for any $C\geq 1$. 

\begin{proposition}{\label{prop1}}
The following results hold for a Banach space $X$.
\begin{enumerate}
\item If every normalized basic sequence in $X$ is Haar-$\ell^+_1$, then every spreading model of $X$ is equivalent to the unit vector basis $\ell_{1}$. 

\item Let $(e_{i})$ be a normalized basic sequence in $X$. If every asymptotic model generated by an array of normalized block bases of $(e_{i})$ is equivalent to the unit vector basis of $\ell_{1}$, then $(e_i)$ is a Haar-$\ell^+_1$ sequence. In particular, $X$ has the Lebesgue property if every asymptotic model of $X$ is equivalent to the unit vector basis of $\ell_1$.
\end{enumerate}
\end{proposition}
\begin{proof}
(1) Let $(e_{i})$ be a normalized and good basic sequence in $X$ that generates the spreading model $(v_{i})$. Consider the normalized difference sequence $(x_i)=\Big(\frac{e_{2i}-e_{2i+1}}{\|e_{2i}-e_{2i+1}\|}\Big)$. Then, $(x_i)$ generates the spreading model which is equivalent to the difference sequence $(w_i)=\Big(\frac{v_{2i}-v_{2i+1}}{\|v_{2i}-v_{2i+1}\|}\Big)$ and it follows from a well-known argument that the sequence $(w_i)$ is $K_{u}$-unconditonal for some $K_{u}\geq 1$. It is therefore enough to show that the sequence $(w_i)$ is equivalent to the unit vector basis of $\ell_1$. Fix a Haar system of partitions $(A^{n}_{j})_{j=0,n\in\N}^{2^{n}-1}$ of $\N$ and note that, because $(x_i)$ is Haar-$\ell^+_1$, there exists $C\ge 1$ such that for all $n\in\N$ there is $m\geq n$ so that 
\bes \frac{1}{2^{m}}\left\Vert\sum_{j=0}^{2^{m}-1}x_{i_{j}} \right\Vert \ge \frac{1}{2C} \ees
for some $2^m\le i_j\in A^m_j$, $0\le j\le 2^m-1$. Now, because $(x_i)$ is a good sequence and $2^m\le i_0<\ldots < i_{2^{m}-1}$, it follows that $(x_{i_{j}})_{j=0}^{2^{m}-1}$ is $2$-equivalent to $(w_{i})_{i=1}^{2^{m}}$ so that
\bes \frac{1}{2C}\leq\frac{1}{2^{m}}\left\Vert\sum_{j=0}^{2^{m}-1}x_{i_{j}}\right\Vert \leq \frac{2}{2^{m}}\left\Vert \sum_{i=1}^{2^{m}}w_{i}\right\Vert \ees
and thus $\frac{1}{4C}\leq \frac{1}{2^{m}}\left\Vert\sum_{i=1}^{2^{m}}w_{i}\right\Vert$. Again by a well-known argument, it now follows that $(w_i)$ is equivalent to the unit vector basis of $\ell_1$ as is required. To see this, let $(a_i)_{i=1}^{2^m}$ be non-negative scalars and, by averaging the vectors $\sum_{i=1}^{2^m}a_{\rho(i)}w_i$ where $\rho$ is a cyclic permutation of $\{1, \ldots, 2^m\}$, we get that

\bes
\left\Vert\sum_{i=1}^{2^{m}}a_iw_{i} \right\Vert \ge
\frac{1}{2}\frac{\sum_{j=1}^{2^m}a_j}{2^m}\left\Vert\sum_{i=1}^{2^{m}}w_{i} \right\Vert 
\ge \frac{1}{8C}\sum_{j=1}^{2^m}a_j.
\ees 
and it then follows from $K_{u}$-unconditionality that for arbitrary $(a_i)_{i=1}^{2^m}$ there is
\bes
\left\Vert\sum_{i=1}^{2^{m}}a_iw_{i} \right\Vert \ge
\frac{1}{8CK_{u}}\sum_{i=1}^{2^m}|a_i|
\ees 
so that the proof of (1) is complete.

(2) In order to prove the second claim, let $(e_{i})$ be a normalized basic sequence in $X$ and assume that every asymptotic model of $X$ with respect to an array of normalized block bases of $(e_{i})$ is equivalent to $\ell_{1}$. Fix a Haar system of partitions $(A_{j}^{n})_{j=0,n\in\N}^{2^{n}-1}$ of $\N$ and let $(A_k)_{k=1}^{\infty}$ be the natural enumeration of $(A_{j}^{n})_{j=0,n\in\N}^{2^{n}-1}$.  Consider the array 
$(e^k_j)$ where each row $(e^k_j)_{j=1}^{\infty}$ is the subsequence $(e_i)_{i\in A_k}$. 
It follows from the assumption that, by passing to a subarray if necessary, $(e^k_j)$ generates an asymptotic model that is $C$-equivalent to the unit vector basis of $\ell_1$ for some $C\geq 1$. This trivially implies that, for all $n$, there is 
\bes
\frac{1}{2^{n}}\left\Vert\sum_{j=1}^{2^{n}}e_{i_j} \right\Vert \ge
\frac{1}{2C}
\ees
for some $i_j\in A^n_{j}$, and this completes the proof of (2).
\end{proof}

\noindent{\bf Remarks.}
 
(1) It is worth noting that the Haar-$\ell^+_1$ condition for all basic sequences does not guarantee that every spreading model is uniformly equivalent to the unit vector basis of $\ell_{1}$. To see this consider the following example. Let $(e_i)$ and $(f_i)$ be the unit vector bases of $\ell_1$ and of $\ell_2$, respectively. For each $k\in\N$, let $X_k\subseteq (\ell_1\oplus \ell_2)_{\ell_1}$ spanned by $(\frac{1}{k}e_i+f_i)_i$. Clearly, $X_k$ is isomorphic to $\ell_1$. In fact, $(\frac{1}{k}e_i+f_i)$ is $k$-equivalent to $(e_i)$, and $k$ is the best constant. Let $X=\Big(\sum_{k=1}^{\infty}X_k\Big)_{\ell_1}$. Then $X$ has the Lebesgue property since the property is stable under $\ell_{1}$-sums (see \cite{dRF}), and by construction, there is no uniform constant for $\ell_1$-spreading models of $X$. Note also that $X$ embeds into $L_1$.

(2) The converse to the first statement in Proposition \ref{prop1} is false. In particular, both Haydon in \cite{Haydon} and Naralenkov in \cite{Naralenkov} construct examples of Schur spaces (hence all spreading models are $\ell_1$) that fail the Lebesgue property.

(3) The converse to the second statement in Proposition \ref{prop1} is also false. Indeed, consider the space $T(T)=(T\oplus T\oplus \ldots)_T$ where $T$ is the Tsirelson space. Let $(e^i_j)_{i,j}$ be the unit vector basis of $T(T)$ in the natural order where $(e^i_j)_j$ is the unit vector basis of $T$ in the $i$th component of the sum. If we consider $(e^i_j)_{i,j}$ as an infinite array where $i$th row is $(e^i_j)_j$, then clearly for all $k\le n_1<\ldots<n_k$, the diagonal sequence $(e^i_{n_i})_{i=1}^k$ is 1-equivalent to the unit vector basis of $T$. That is, the array generates the unit vector basis of $T$ as an asymptotic model. On the other hand, the following proposition shows that $T(T)$ has the Lebesgue property. The proof of this proposition is similar to the argument that is used in \cite{dRF} to show that the Lebesgue property is stable under $\ell_{1}$ sums and to the argument that is used in \cite{Gaebler} to show that a Banach space has the Lebesgue property if it is asymptotic-$\ell_{1}$ in the coordinate-free sense of \cite{MMT}.

\begin{proposition}{\label{Tsirelson_sum}}
Let $\{X_{i}\}_{i=1}^{\infty}$ be a collection of Banach spaces such that $X_{i}$ has the Lebesgue property for each $i\in\N$. Then, $\left(\sum_{i=1}^{\infty}X_{i}\right)_{T}$ has the Lebesgue property.
\end{proposition}
\begin{proof}
Let $f:[0,1]\to \left(\sum_{i=1}^{\infty}X_{i}\right)_{T}$ be a bounded function that is discontinuous on a set $H$ of positive Lebesgue measure. We will prove that $f$ cannot be Riemann-integrable. Note first that $H=\bigcup_{n=1}^{\infty}H_{n}$ where
\bes H_{n}=\left\{t\in[0,1] \;\middle\vert\; \inf_{\delta>0}\sup\{\|f(s)-f(s')\| \mid s,s'\in[t-\delta,t+\delta]\cap[0,1]\}\geq\frac{1}{2}\right\} \ees
so that $\mu(H_{n_{0}})>0$ for some $n_{0}\in\N$, else $\mu(H)=0$. Next, let $P_{i}$ be the continuous linear projection of $\left(\sum_{i=1}^{\infty}X_{i}\right)_{T}$ onto $X_{i}$ and we may assume without loss of generality that
\bes G_{i}=\{t\in[0,1] \mid (P_{i}\circ f):[0,1]\to X_{i} \text{ is discontinuous at } t\} \ees
has Lebesgue measure zero for every $i\in\N$. Indeed, if $\mu(G_{i})>0$ for some $i\in\N$, then $P_{i}\circ f$ cannot be Riemann-integrable because $X_{i}$ has the Lebesgue property and this implies that $f$ cannot be Riemann-integrable so that there is nothing more to prove. Thus, $G=\bigcup_{i=1}^{\infty}G_{i}$ also has Lebesgue measure zero. Now, let $P=\left(\frac{k}{2^{N}}\right)_{k=0}^{N}:=(p_{k})_{k=0}^{N}$ be an arbitrary dyadic partition of $[0,1]$ and define the set
\bes A=\{ k \mid \mu[(p_{k-1},p_{k})\cap (H_{n_{0}}\setminus G)]>0\} \ees
which is non-empty because $H_{n_{0}}\setminus G$ has positive Lebesgue measure and let $k_{1},\ldots,k_{r}$ be the members of $A$. 

Let $r=m_{0}\in\N$ and let $\ep>0$ be given. Then, let $s_{1}\in(p_{k_{1}-1},p_{k_{1}})\cap (H_{n_{0}}\setminus G)$ and choose $t_{1},t'_{1}\in(p_{k_{1}-1},p_{k_{1}})$ so that $\|z_{1}\|:=\|f(t_{1})-f(t'_{1})\|\geq\frac{1}{2_{n_{0}}}$ and so that $\sum_{i=1}^{m_{0}}\|(P_{i}\circ f)(t_{1})-(P_{i}\circ f)(t'_{1})\|<\ep$ because $P_{i}\circ f$ is continuous at $s_{1}$ for each $i\in\N$. Choose $m_{1}>m_{0}$ so that $\sum_{i=m_{1}}^{\infty}\|P_{i}(z_{1})\|<\ep$.

Now, let $s_{2}\in(p_{k_{2}-1},p_{k_{2}})\cap (H_{n_{0}}\setminus G)$ and choose $t_{2},t'_{2}\in(p_{k_{2}-1},p_{k_{2}})$ so that $\|z_{2}\|:=\|f(t_{2})-f(t'_{2})\|\geq\frac{1}{2n_{0}}$ and so that $\sum_{i=1}^{m_{1}}\|(P_{i}\circ f)(t_{2})-(P_{i}\circ f)(t'_{2})\|<\frac{\ep}{2}$ because $P_{i}\circ f$ is continuous at $s_{2}$ for each $i\in\N$. Choose $m_{2}>m_{1}$ so that $\sum_{i=m_{2}}^{\infty}\|P_{i}(z_{2})\|<\frac{\ep}{2}$.

Continue in this manner for the remaining $k_{3},\ldots,k_{r}$ sub-intervals that correspond to the members of $A$ to find points $t_{l},t'_{l}\in(p_{k_{l}-1},p_{k_{l}})$ so that $\|z_{l}\|:=\|f(t_{l})-f(t'_{l})\|\geq\frac{1}{2n_{0}}$ and positve integers $m_{l-1}<m_{l}$ so that $\sum_{i=1}^{m_{l-1}}\|(P_{i}\circ f)(t_{l})-(P_{i}\circ f)(t'_{l})\|<\frac{\ep}{2^{l-1}}$ and $\sum_{i=m_{l}}^{\infty}\|P_{i}(z_{l})\|<\frac{\ep}{2^{l-1}}$. Then, choose interior tags $t_{k}=t'_{k}$ for each $k\notin A$ and note that
\begin{align*} \frac{1}{2^{N}}\left\Vert \sum_{l=1}^{r}z_{l}\right\Vert &\geq \frac{1}{2^{N}}\left\Vert\sum_{l=1}^{r}\sum_{i=m_{l-1}+1}^{m_{l}-1}P_{i}(z_{l})\right\Vert -4\ep \\ &\geq \frac{1}{2^{N}}\frac{1}{2}\sum_{l=1}^{r}\left\Vert\sum_{i=m_{l-1}+1}^{m_{l}-1}P_{i}(z_{l})\right\Vert-4\ep \\ &\geq\frac{1}{2^{N}}\frac{1}{2}\sum_{l=1}^{r}\left(\|z_{l}\|-4\ep\right)-4\ep\geq \frac{1}{2}\left(\frac{1}{2n_{0}}-4\ep\right)\mu(H_{n_{0}}\setminus G)-4\ep \end{align*}
by the definition of the Tsirelson norm. This completes the proof because we are free to take $\ep>0$ as small as we want and, thus, $f$ cannot be Riemann-integrable.
\end{proof}

(4) In general, one cannot improve the asymptotic structure of a space with the Lebesgue property by passing to subspaces as the next example shows.

\subsection{Argyros-Motakis space $X_{iw}$} The Banach space $X_{iw}$ that is constructed by Argyros and Motakis in \cite{AM} has the remarkable property that, {\em while all spreading models are uniformly equivalent to the unit vector basis of $\ell_1$, every subspace admits an asymptotic model that is generated by an array of normalized block bases that is equivalent to the unit vector basis of $c_0$.} This extremal separation between asymptotic structures extends to the separation between the Lebesgue property and asymptotic models as well. Indeed, we show below that the space $X_{iw}$ actually has the Lebesgue property despite the presence of $c_{0}$ asymptotic models in every subspace. The proof is a slight generalization of the proof of \cite[Proposition 4.1]{AM}, however, we choose to include the details here for completeness because there are several additional considerations.

Recall that the Schreier sets form an increasing sequence of collections of finite subsets of $\N$ that are defined inductively as follows:
\bes S_{0}=\{\{i\} \mid i\in\N\} \text{ and } S_{1}=\{ F\subset\N \mid |F|\leq\min(F)\} \ees
where $|F|$ is the cardinality of $F\subset\N$ and for all $n$, 
\bes S_{n+1}=\left\{ F\subset\N \;\middle\vert\; F=\bigcup_{i=1}^{d}F_{i}, F_{1}<\ldots<F_{d},\  { F_{i}\in S_{n},\   d\leq\min(F_{1})} \right\}.\ees
A sequence $x_{1}<\ldots<x_{r}$ of vectors in $c_{00}$ is said to be $S_{n}$-admissible if $\{\ms(x_{i}) \mid 1\leq i\leq r\}\in S_{n}$. Let $(a_{k})_{k}$ be a strictly increasing sequence of positive integers with $a_{1}=2$ and let $(b_{k})_{k}$ be another strictly increasing sequence of positive integers such that $b_{1}=1$ and such that $b_{k+1}>b_{k}\log(a^{2}_{k+1})$ (this technical condition is given in \cite{AM} but it is not required for our proof below). We now construct the norming functionals for $X_{iw}$. 

\begin{definition}{\label{Xiw_conventions}}
The norming set $G\subset c_{00}$ is to satisfy the following requirements.
\begin{enumerate}
\item Let $k_{1},\ldots,k_{l}$ be positive integers and let $f_{1}<\ldots<f_{r}\in G$ be $S_{b_{k_{1}}+\cdots+b_{k_{l}}}$-admissible functionals. Then,
\bes f=\frac{1}{a_{k_{1}}\cdots a_{k_{l}}}\sum_{q=1}^{r}f_{q} \in G \ees
and it is said to have weight $w(f)=a_{k_{1}}\cdots a_{k_{l}}$. Also, the functionals $\pm e_{i}^{*}$ are defined to have infinite weight.
\item A sequence $f_{1}<f_{2}<\ldots$ of weighted functionals in $G$ is said to be very fast growing if $\Ms(f_{q-1})<w(f_{q})$ for $q>1$.
\end{enumerate}
\end{definition}

\noindent The desired norming set $G=W_{iw}$ for $X_{iw}$ is then the smallest subset of $c_{00}$ such that $\pm e_{i}^{*}\in W_{iw}$ and such that $f=\frac{1}{a_{k_{1}}\cdots a_{k_{l}}}\sum_{q=1}^{r}f_{q}\in W_{iw}$ if the weighted functionals $(f_{q})_{q=1}^{r}$ are in $W_{iw}$, are $S_{b_{k_{1}}+\cdots+b_{k_{l}}}$-admissible, and are very fast growing. The Banach space $X_{iw}$ is the completion of $c_{00}$ with respect to the norm $\|x\|=\sup\{f(x) \mid f\in W_{iw}\}$.

\begin{proposition}{\label{prop2}}
The Banach space $X_{iw}$ has the Lebesgue property.
\end{proposition}
\begin{proof}
Suppose that $X_{iw}$ does not have the Lebesgue property. Then, there exists from Proposition \ref{Haydon-Odell} a normalized block basis $(x_{j})$ such that the function $g:[0,1]\to X_{iw}$ that is defined by
\bes g(t)=\begin{cases} x_{j} &\mbox{if } t=d_{j}\in D \\ 0 &\mbox{if } t\notin D \end{cases} \ees
is Riemann-integrable, where $D$ is as usual the set of dyadic rational numbers in $(0,1)$. For every $j\in\N$, choose a weighted functional $f_{j}\in W_{iw}$ such that $f_{j}(x_{j})=1$ and such that $\text{supp}(f_{j})\subset\text{supp}(x_{j})$. There are now two cases.

\vspace{0.25cm}

\noindent\textbf{Case 1:} $\sup_{j\in\N}w(f_{j})=L<\infty$

\vspace{0.25cm}

\noindent In this case, for each $s\in \Lambda:=[0,1]\setminus D$, fix a sequence $(d^{s}_{j})$ of dyadic rational numbers that converges to $s$. Then, by passing to a subsequence if necessary, we may assume that there is some $n\in\N\cap[1,L]$ such that $w(f^{s}_{j})=n$ for each $j\in\N$. We therefore define the set $\Lambda_{n}=\{ s \mid w(f_{j}^{s}) = n\}$ for each $1\leq n\leq L$ and note that
\bes 1=\mu(\Lambda)=\mu^{*}(\Lambda)\leq\sum_{n=1}^{L}\mu^{*}(\Lambda_{n}) \ees
because $\Lambda=\cup_{n=1}^{L}\Lambda_{n}$, where $\mu^{*}$ denotes the Lebesgue outer measure. Evidently, it follows that $\mu^{*}(\Lambda_{n_{0}})>0$ for some $n_{0}\in\N\cap[1,L]$. 

Let $\ep>0$ and choose $n\in\N$ so that $\frac{1}{n}\left\Vert\sum_{i=1}^{n}g(t_{i})\right\Vert<\ep$ for every choice $t_{i}\in\left(\frac{i-1}{n},\frac{i}{n}\right)$ of interior tags. Then, write $\left\{\left(\frac{i-1}{n},\frac{i}{n}\right)\cap \Lambda_{n_{0}}\neq\emptyset\right\}=\{i_{1},\ldots,i_{r}\}$ where this set is non-empty because $\mu^{*}(\Lambda_{n_{0}})>0$ and note that $\frac{r}{n}\geq\mu^{*}(\Lambda_{n_{0}}\setminus\{\frac{i}{n}\}_{i=0}^{n})=\mu^{*}(\Lambda_{n_{0}})$ because $\Lambda_{n_{0}}\setminus\{\frac{i}{n}\}_{i=0}^{n}\subset\cup_{l=1}^{r}\left(\frac{i_{l}-1}{n},\frac{i_{l}}{n}\right)$. Now, choose interior tags $t_{1}=d^{s_{1}}_{j}\in\left(\frac{i_{1}-1}{n},\frac{i_{1}}{n}\right),\ldots, t_{r}=d^{s_{r}}_{j}\in\left(\frac{i_{r}-1}{n},\frac{i_{r}}{n}\right)$ where $s_{l}\in\left(\frac{i_{l}-1}{n},\frac{i_{l}}{n}\right)\cap\Lambda_{n_{0}}$ for $1\leq l\leq r$ so that $x_{1}<\ldots<x_{r}$ are normalized blocks with $2^{n}\leq\min\supp(x_{1})$.

Let $f_{l}\in W_{iw}$ be the norming functionals that correspond to these block vectors and note that
\bes f_{l}=\frac{1}{a_{k_{1}}\cdots a_{k_{p}}}\sum_{q=1}^{d_{l}}f^{l}_{q} \ees
where $a_{k_{1}}\cdots a_{k_{p}}=n_{0}$, and where $f^{l}_{1}<\ldots<f^{r}_{d_{l}}$ are very fast growing and $S_{b_{k_{1}}+\cdots+b_{k_{p}}}$-admissible. Note that $|f^{l}_{q}(x_{l})|\leq 1$ because $f^{l}_{q}$ is in the unit ball of $X_{iw}^{*}$ and, in particular, there is
\bes \frac{1}{a_{k_{1}}\cdots a_{k_{p}}}\sum_{q=2}^{d_{l}}f^{l}_{q}(x_{l})=f_{l}(x_{l})-\frac{1}{a_{k_{1}}\cdots a_{k_{p}}}f^{l}_{1}(x_{l}) \geq 1-\frac{1}{2}=\frac{1}{2} \ees
and, what is more, it follows that the functionals
\bes f^{1}_{2}<\ldots<f^{1}_{d_{1}}<f^{2}_{2}<\ldots<f^{2}_{d_{2}}<\ldots<f^{r}_{2}<\ldots<f^{r}_{d_{r}} \ees
are very fast growing and $S_{b_{k_{1}}+\cdots+b_{k_{p}}+1}$-admissible. This now implies that
\bes F=\frac{1}{a_{k_{1}}\cdots a_{k_{p}}a_{1}}\sum_{l=1}^{r}\sum_{q=2}^{d_{l}}f^{l}_{q} \in W_{iw} \ees
and, consequently, there is the estimate
\begin{align*} \frac{1}{n}\left\Vert\sum_{i=1}^{n}g(t_{i})\right\Vert & =\frac{1}{n}\left\Vert\sum_{l=1}^{r}g(t_{l})\right\Vert =\frac{1}{n}\left\Vert\sum_{l=1}^{r}x_{l}\right\Vert\geq\frac{1}{n}F\left(\sum_{l=1}^{r}x_{l}\right) \\ &=\frac{1}{n}\frac{1}{a_{k_{1}}\cdots a_{k_{p}}a_{1}}\sum_{v=1}^{r}\sum_{q=2}^{d_{v}}f^{v}_{q}\left(\sum_{l=1}^{r}x_{l}\right) \\ &=\frac{1}{a_{1}n}\sum_{l=1}^{r}\frac{1}{a_{k_{1}}\cdots a_{k_{p}}}\sum_{q=2}^{d_{l}}f^{l}_{q}(x_{l})\geq \frac{r}{4n}\geq\frac{\mu^{*}(\Lambda_{n_{0}})}{4} \end{align*}
where $t_{i}\notin D$ for $i\notin\{i_{1},\ldots,i_{r}\}$ and this now contradicts the fact that $g$ is Riemann-integrable for $\ep$ sufficiently small so that the argument for Case 1 is complete. 

\vspace{0.25cm}

\noindent\textbf{Case 2:} $\sup_{j\in\N} w(f_{j})=\infty$

\vspace{0.25cm}

\noindent Note first that if $\sup\{w(f_{j}) \mid d_{j}\in I\}<\infty$ for some dyadic sub-interval $I\subset[0,1]$, then the argument for Case 1 again yields a contradiction to the Riemann integrability of $h(t)=\mathbbm{1}_{I}(t)g(t)$. We may therefore assume that $\sup\{w(f_{j}) \mid d_{j}\in I\}=\infty$ for every dyadic sub-interval $I\subset[0,1]$. Now, given $n\in\N$, choose $\tilde{d}_{1},\ldots,\tilde{d}_{2^{n}}$ in the interiors of the $2^{n}$ dyadic sub-intervals of $[0,1]$ so that 
\bes \tilde{x}_{1}<\ldots<\tilde{x}_{2^{n}} \qquad \tilde{x}_{i}=g(\tilde{d}_{i}) \ees
is $S_{1}$-admissible and so that the norming functionals $\tilde{f}_{1}<\ldots<\tilde{f}_{2^{n}}$ are very fast growing. It follows that
\bes F=\frac{1}{a_{1}}\sum_{i=1}^{2^{n}}\tilde{f}_{i}\in W_{iw} \ees
and, consequently, there is the estimate
\begin{align*} \frac{1}{2^{n}}\left\Vert\sum_{i=1}^{2^{n}}g(\tilde{d}_{i})\right\Vert &= \frac{1}{2^{n}}\left\Vert\sum_{i=1}^{2^{n}}\tilde{x}_{i}\right\Vert \geq F\left(\sum_{i=1}^{2^{n}}\tilde{x}_{i}\right) \\ &=\frac{1}{2^{n}}\frac{1}{a_{1}}\sum_{i=1}^{2^{n}} \tilde{f}_{i}\left(\sum_{i'=1}^{2^{n}}\tilde{x}_{i'}\right) = \frac{1}{2^{n}}\frac{1}{a_{1}}\sum_{i=1}^{2^{n}}\tilde{f}_{i}(\tilde{x}_{i}) = \frac{1}{2}\end{align*}
and this now contradicts the fact that $g$ is Riemann-integrable for $n$ sufficiently large so that the proof is complete.
\end{proof}

\section{The $L_{1}[0,1]$ Theorem of Pelczynski and da Rocha Filho}{\label{Sec5}}

In this section, we give an alternate and direct proof of the following theorem of Pelczynski and da Rocha Filho.

\begin{theorem}{\label{MR2}}
A subspace $X\subset L_{1}[0,1]$ has the Lebesgue property if every spreading model of $X$ is equivalent to $\ell_{1}$.
\end{theorem}

Our proof uses a fundamental result that is due to Dor in \cite{Dor}. This result is stated below as a lemma and it asserts that the sequences $(f_i)\subset L_1[0,1]$ that are equivalent to the unit vector basis of $\ell_1$ are essentially disjointly supported.   

\begin{lemma}{\label{Dor}}
Let $n\in\N$, let $f_{1},\ldots,f_{n}$ be functions in the unit sphere of $L_{1}[0,1]$, and let $\theta\in(0,1]$. If, for every choice $(a_{i})_{i=1}^{n}$ of scalars, there is the estimate
\bes \theta\sum_{i=1}^{n}|a_{i}|\leq\left\Vert\sum_{i=1}^{n}a_{i}f_{i}\right\Vert \ees
then there exist disjoint and Lebesgue measurable subsets $A_{1},\ldots,A_{n}$ of $[0,1]$ such that $\|f_{i}|_{A_{i}}\|=\int_{A_{i}}|f_{i}(t)|dt\geq\theta^{2}$ for each $1\leq i\leq n$.
\end{lemma}

We now give our proof of Theorem \ref{MR2}. The beginning of this proof is similar to the first case of Proposition \ref{prop2}.

\begin{proof}[Proof of Theorem \ref{MR2}]
Suppose for a contradiction that $X$ does not have the Lebesgue property. It follows from from Proposition \ref{Haydon-Odell} that there exists a normalized basic sequence $(x_{j})_{j}$ so that the everywhere discontinuous function $f:[0,1]\to X$ defined by
\bes f(t)=\begin{cases} x_{j} &\mbox{if } t=d_{j}\in D \\ 0 &\mbox{if } t\notin D\end{cases} \ees
is Riemann-integrable, where $D$ is the set of dyadic rational numbers in $(0,1)$. Next, for each $s\in \Lambda:=[0,1]\setminus D$, let $(d^{s}_{j})$ be a sequence of members of $D$ that converges to $s$ and, by passing to a subsequence if necessary, we may assume that $(f(d^{s}_{j}))$ is a good sequence in $X$ that generates a spreading model $(v^{s}_{j})$ which is equivalent to the unit vector basis of $\ell_{1}$. Then, for each $m\in\N$, define the set $\Lambda_{m}=\{s\in \Lambda \mid (v_{j}^{s}) \overset{m}{\sim}\ell_{1}\}$ where $\overset{m}{\sim}$ denotes the equivalence with constant $m$. It follows that $\Lambda=\cup_{m=1}^{\infty}\Lambda_{m}$ so that there is
\bes 1=\mu(\Lambda)=\mu^{*}(\Lambda)\leq\sum_{m=1}^{\infty}\mu^{*}(\Lambda_{m}) \ees
where $\mu^{*}$ denotes the Lebesgue outer measure and, evidently, we have that $\mu^{*}(\Lambda_{m_0})>0$ for some $m_0\in\N$.

Let $\ep>0$. Then, choose $n\in\N$ large enough so that $\frac{1}{n}<\ep$ and so that, because $f$ is Riemann-integrable, there is the estimate
\bes \frac{1}{n}\left\Vert\sum_{i=1}^{n}(f(t_{i})-f(t'_{i}))\right\Vert<\ep \ees
for every choice $t_{i},t'_{i}\in\left(\frac{i-1}{n},\frac{i}{n}\right)$ of interior tags. Note moreover that the function takes zero values in every interval so that actually we have stronger inequality that 
\bes \frac{1}{n}\max_{\ep_i=\pm 1}\left\Vert\sum_{i=1}^{n}\ep_if(t_{i})\right\Vert<\ep \ees
 for every choice $t_{i}\in\left(\frac{i-1}{n},\frac{i}{n}\right)$ of interior tags. 

Now, let $\{ i \mid \left(\frac{i-1}{n},\frac{i}{n}\right)\cap\Lambda_{m_{0}}\neq \emptyset \}=\{i_{1},\ldots,i_{r}\}$ where this set is non-empty because $\mu^{*}(\Lambda_{m_{0}})>0$. It follows that $\frac{r}{n}\ge \mu^*(\Lambda_{m_0}) = \mu^{*}(\Lambda_{m_{0}}\setminus\{\frac{i}{n}\}_{i=0}^{n})$ because $\Lambda_{m_{0}}\setminus\{\frac{i}{n}\}_{i=0}^{n}\subset\cup_{k=1}^{r}\left(\frac{i_{k}-1}{n},\frac{i_{k}}{n}\right)$. Then, for $2\leq k\leq r$, fix a point $s_{k}\in\left(\frac{i_{k}-1}{n},\frac{i_{k}}{n}\right)\cap\Lambda_{m_{0}}$ so that there is an infinite sequence $(d^{s_{k}}_{j})$ that converges to $s_{k}$ where $(f(d^{s_{k}}_{j}))$ is a good sequence in $X$ that generates a spreading model $(v^{s_{k}}_{j})$ which is $m_{0}$-equivalent to the unit vector basis of $\ell_{1}$. 

Fix an arbitrary point $\tau_{1}\in \left(\frac{i_1-1}{n},\frac{i_1}{n}\right)\cap D$. Next, consider the sequence $(d^{s_{2}}_{j})$ that converges to $s_{2}\in\left(\frac{i_{2}-1}{n},\frac{i_{2}}{n}\right)\cap\Lambda_{m_{0}}$ and note that, for every $N\in\N$ and for every choice $(a_{l})_{l=1}^{N}$ of scalars, there is
\bes \frac{1}{2m_{0}}\sum_{l=1}^{N}|a_{l}| \leq \left\Vert\sum_{l=1}^{N}a_{l}f(d^{s_{2}}_{j_{l}})\right\Vert \ees
for $N\leq j_{1}<\ldots<j_{N}$ because $(f(d^{s_{2}}_{j_{l}}))_{l=1}^{N}$ is $2$-equivalent to $(v_{j}^{s_{2}})_{j=1}^{N}$. It therefore follows from Lemma \ref{Dor} that there exist disjoint and Lebesgue measurable sets $U^{2}_{1},U^{2}_{2},\ldots, U^{2}_{N}$ in $[0,1]$ so that $\frac{1}{4m_{0}^{2}}\leq \|f(d^{s_{2}}_{j_{l}})|_{U^{2}_{l}}\|=\int_{U^{2}_{l}}|f(d^{s_{2}}_{j_{l}})(w)|dw$ for each $1\leq l\leq N$. Moreover, we have that
\bes 1=\|f(\tau_{1})\|\geq\sum_{l=1}^{N}\int_{U^{2}_{l}}|f(\tau_{1})(w)|dw \ees
because the sets $U^{2}_{1},U^{2}_{2},\ldots, U^{2}_{N}$ are disjoint and, by choosing $N$ large enough, it follows that $\int_{U_{l}^{2}}|f(\tau_{1})(w)|dw<\frac{\ep}{2}$ for some $1\leq l\leq N$. For such an $l$, set $U_{2}=U^{2}_{l}$ and $\tau_{2}=d^{s_{2}}_{j_{l}}$.

Next, consider the sequence $(d^{s_{3}}_{j_{l}})$ that converges to $s_{3}\in\left(\frac{i_{3}-1}{n},\frac{i_{3}}{n}\right)\cap\Lambda_{m_{0}}$ and note as above that, for every $N\in\N$ and for every choice $(a_{l})_{l=1}^{N}$ of scalars, there is
\bes \frac{1}{2m_{0}}\sum_{l=1}^{N}|a_{l}| \leq \left\Vert\sum_{l=1}^{N}a_{l}f(d^{s_{3}}_{j_{l}})\right\Vert \ees
for $N\leq j_{1}<\ldots<j_{N}$ because $(f(d^{s_{3}}_{j_{l}}))_{l=1}^{N}$ is $2$-equivalent to $(v_{j}^{s_{3}})_{j=1}^{N}$. It therefore follows from Lemma \ref{Dor} that there exist disjoint and Lebesgue measurable sets $U^{3}_{1},U^{3}_{2},\ldots, U^{3}_{N}$ in $[0,1]$ so that $\frac{1}{4m_{0}^{2}}\leq \|f(d^{s_{3}}_{j_{l}})|_{U^{3}_{l}}\|=\int_{U^{3}_{l}}|f(d^{s_{3}}_{j_{l}})(w)|dw$ for each $1\leq l\leq N$. Moreover, we have that
\begin{align*} 1&=\|f(\tau_{1})\|\geq\sum_{l=1}^{N}\int_{U^{3}_{l}}|f(\tau_{1})(w)|dw, \\ 1&=\|f(\tau_{2})\|\geq\sum_{l=1}^{N}\int_{U^{3}_{l}}|f(\tau_{2})(w)|dw \end{align*}
because the sets $U^{3}_{1},U^{3}_{2},\ldots, U^{3}_{N}$ are disjoint and, by choosing $N$ large enough, it follows that $\int_{U_{l}^{3}}|f(\tau_{1})(w)|dw<\frac{\ep}{4}$ and $\int_{U_{l}^{3}}|f(\tau_{2})(w)|dw<\frac{\ep}{4}$ for some $1\leq l\leq N$. For such an $l$, set $U_{3}=U^{3}_{l}$ and $\tau_{3}=d^{s_{3}}_{j_{l}}$.

We then continue in this manner to find Lebesgue measurable subsets $U_{k}$ and interior tags $\tau_{k}=d^{s_{k}}_{j_{l}}$ as above so that
\begin{align*}\int_{U_{q}}|f(\tau_{k})(w)|dw&<\frac{\ep}{2^{q-1}} \text{ for } k< q, \text{ whereas } \int_{U_{k}}|f(\tau_{k})(w)|dw\geq\frac{1}{4m_{0}^{2}}
\end{align*}
for each $2\leq k\leq r$. Then, by the square function estimate (Khintchine) (cf. \cite[Theorem 6.2.13]{AK}), we have that
\begin{align*} \ep&>\max_{\ep_i=\pm}\frac{1}{n}\left\Vert\sum_{i=1}^{n}\ep_i f(\tau_{i})\right\Vert\geq \underset{\ep_i=\pm}{\text{avg}}\frac{1}{n}\left\Vert\sum_{i=1}^{n}\ep_i f(\tau_{i})\right\Vert \\ &\geq\frac{1}{n}\left\Vert\left(\sum_{i=1}^{n}|f(\tau_{i})|^{2}\right)^{\frac{1}{2}}\right\Vert=\frac{1}{n}\left\Vert\left(\sum_{k=1}^{r}|f(\tau_{k})|^{2} \right)^{\frac{1}{2}}\right\Vert\end{align*}
where we take $\tau_i\not\in D$ for $i\notin\{i_{1},\ldots,i_{r}\}$. Now, write $V_{k}=U_{k}\setminus\left(\cup_{q>k}U_{q}\right)$ for $2\leq k \leq r$ so that $V_{2},\ldots,V_{r}$ are pairwise disjoint and the above quantity is bounded below by
\begin{align*} &\frac{1}{n}\sum_{z=2}^{r}\int_{V_{z}}\left(\sum_{k=1}^{r}|f(\tau_{k})|^{2} \right)^{\frac{1}{2}}(w)dw \geq \frac{1}{n}\sum_{k=2}^{r}\int_{V_{k}}|f(\tau_{k})(w)|dw \\ &\geq\frac{1}{n}\sum_{k=2}^{r}\left(\int_{U_{k}}|f(\tau_{k})(w)|dw-\sum_{q>k}\int_{U_{q}}|f(\tau_{k})(w)|dw\right) \\ &\geq\frac{1}{n}\sum_{k=2}^{r}\left(\frac{1}{4m_{0}^{2}}-\ep\right)\geq \left(\frac{1}{4m_{0}^{2}}-\ep\right)\left(\frac{r}{n}-\frac{1}{n}\right)\geq\left(\frac{1}{4m_{0}^{2}}-\ep\right)\left(\mu^{*}(\Lambda_{m_{0}})-\ep\right)\end{align*}
and this is a contradiction for small enough $\ep>0$ so that the proof is complete.
\end{proof}

%%%%%%%%%%%%%%%%%%%%%%%%%%%%%%%%%%%%%%%%%%%%%%%%%%%%%%%%%%%%%%%%%%%%%%

\end{document}